\documentclass[9pt,twocolumn,twoside]{ksarticle}

\usepackage{bm}
\usepackage{amssymb}
\usepackage{xcolor}
\usepackage[normalem]{ulem}
\usepackage{booktabs}
\usepackage{comment}
\usepackage{hyperref}
\usepackage{authblk}
\usepackage{microtype}
\usepackage[switch]{lineno}

\DeclareRobustCommand{\revdel}[1]{{\let\oldref\ref\def\ref##1{\mbox{\oldref{##1}}}\let\oldautoref\autoref\def\autoref##1{\mbox{\oldautoref{##1}}}\textcolor{red}{\sout{#1}}}}

\newcommand{\cEi}{\mathcal E_i(t)}

\newcommand{\cKik}{\mathcal K_i^k(t)}

\AtBeginDocument{%
  }
\newcounter{subfig}[figure]

\newcommand{\sublabel}[1]{\refstepcounter{subfig}\label{#1}\textbf{(\alph{subfig})}}

\newcounter{sisec}
\renewcommand{\thesisec}{S\arabic{sisec}}

\newcommand{\sisection}[1]{\refstepcounter{sisec}\section*{\thesisec. #1}}
\additionalelement{}
\begin{document}

\title{Mathematics of Herding}

\author[a,1]{Sayomi Kamimoto}
\affil[a]{Department of Mathematics, Howard University, Washington, DC 20059}
\leadauthor{Kamimoto}

\authorcontributions{S.K. designed research, performed research, analyzed data, and wrote the paper, with assistance from Claude (Anthropic) in editing, LaTeX preparation, coding, and algebraic and numerical checks.}
\authordeclaration{The author declares no competing interests.}
\correspondingauthor{\textsuperscript{1}To whom correspondence should be addressed. E-mail: sayomi.kamimoto@howard.edu}

\begin{abstract}
A herd of animals often follows a leader who sets the direction of travel. The leader sometimes turns unexpectedly and the followers reshuffle. Yet the herd regroups itself and stays coherent. Our work aims to understand the herding behavior mathematically. We construct a minimal discrete-time nonequilibrium dynamical system with one-way coupling, and the minimal structure allows an exact analysis. Take the location of each follower in the frame of the leader. The recorded locations form a distribution. The limit of the distribution is the asymptotic spatial measure, whose geometric representation we study. The dynamics of the herd have an exact operator representation with two factors, the operator of a fixed leader-follower pair and the operator of the interaction rules, which reshuffle the pairs. We show that the influence of the leader passes through some follower more than once. As a result, two followers can have reciprocal influence, despite the one-way coupling.
\end{abstract}

\dates{This manuscript was compiled on \today}
\articletype{Research Article}

\maketitle
\thispagestyle{firststyle}
\ifthenelse{\boolean{shortarticle}}{\ifthenelse{\boolean{singlecolumn}}{\abscontentformatted}{\abscontent}}{}

\section*{Introduction}
\label{sec:Intro}
There are group behaviors in nature and engineered systems in which a group of many individuals appears as a single moving object, a herd. The herd has a permanent leader. The leader might direct the herd in a direction the members are not aware of. As the herd moves, each individual other than the leader is at any moment following someone nearby and ahead of it, possibly the leader. Each individual might have followers at any point of time. The herd is then a collection of pairs that is reshuffled continually. 

We translate the group motions of herding into a deterministic discrete-time nonequilibrium dynamical system with one-way coupling. In our model, the herd consists of the leader and \(N\) followers. Each follower follows one individual, and tracing the chain of following from any follower ends at the leader. The pairs at a given time thus form a tree with the leader at the root\cite{BangJensen2009}, and the tree changes as the pair formation changes. The present work shows the existence of a herding system that can be analyzed exactly while producing features that motivate further development.

A tree whose connections coevolve with the positions of the members is an adaptive network~\cite{Gross2008}. The dynamics of adaptive networks are notoriously difficult to analyze~\cite{Sayama2013,Berner2023}. The usual methods rely on an approximation. Mean-field approximation replaces the state of the individuals by their average~\cite{Gross2006,Chen2016}, and separation of time scales imposes the static network while the individuals move~\cite{Kuehn2012,Kuehn2015,Berner2023}. Direct simulation keeps both the individuals and the changing network but does not give an explicit expression of the solution~\cite{Holme2006,Vazquez2008}. We construct, instead, an adaptive network that changes on the same time scale as the motions of the herd. Our method makes no approximation and gives the solution in explicit form.

The present system is out of equilibrium, and the motion of the herd is unbounded. Typically, the analysis of a nonequilibrium system looks for a steady state in a statistical sense, the fraction of time spent in each region. In a case of  collective motion, the steady state is sought in a frame that moves with the group, the center of mass, either as a formation to which the group converges~\cite{Cucker2007,Olfati2006} or as averaged quantities such as order parameters and correlation functions~\cite{Vicsek1995,Ballerini2008}. In contrast, we choose the frame of the leader, and the steady state we consider is the distribution of the relative positions itself.

We record the position of every follower at each time step relative to the leader in coordinates aligned with the direction of travel of the leader. Over time, the recorded positions accumulate into a distribution, and the limit of the distribution is called the asymptotic spatial measure~\cite{Ott2002,Eckmann1985}. The location, extent, and symmetry of the measure in the plane are referred to as the geometric representation of the measure, or simply its geometry or representation, the object of our study.

The dynamics, when written as operators, factor into the operator of a single pair and that of the interaction rules. The measure for the herd is then built from the pair along the changing trees. The influence of the leader is the effect of a past position of the leader on the present position of a follower. The path of the influence is the route along which the effect travels through the pairs. We quantify the influence by the number of steps on the path, one step for each pair the effect travels through. We demonstrate that the influence of the leader reaches a follower through more followers than the herd contains. As a result, the influence between two followers can run in both directions, despite one-way coupling.

\section*{Biological motivation}
\label{sec:motivation}

One-way coupling is the extreme case of nonreciprocal interaction~\cite{Nagy2010,Fruchart2021,Sattari2022,Ivlev2015,Saha2020,You2020}. Dynamical systems of collective motion are typically built on the assumption of reciprocal interactions~\cite{Vicsek1995,Couzin2005,DOrsogna2006,Cucker2007,Vicsek2012}. In models of nonreciprocal active matter, the asymmetry acts between species, and the influence remains in both directions~\cite{Ivlev2015,Fruchart2021,You2020,Saha2020}. One-way coupling between individuals has been studied in consensus models, where the members converge to the state of the leader~\cite{Jadbabaie2003,Olfati2004,Moreau2005}. One-way coupling between individuals without convergence, the class of the present system, has received little attention. 

We turn to the literature in evolutionary biology on tandem running, a leader-follower interaction in termites. For flocks of birds and schools of fish, the following behavior is assumed to be based mainly on vision~\cite{Pearce2014,Strandburg2013}. Tandem running instead requires physical contact: the follower touches the tip of the abdomen of the leader with its palps and antennae~\cite{Mizumoto2024,Mizumoto2025}. Through the contact, information passes between the two individuals. The amount and direction of the transfer are quantified by a conditional mutual information between their time series, known as transfer entropy~\cite{Schreiber2000}. In blind termites, the leader sets the direction and the speed, and information passes almost exclusively from leader to follower~\cite{Valentini2020,Mizumoto2021}. The findings thus justify the one-way coupling of the present system.

The interactions that coordinate collective movement are often pairwise and local, even in large groups~\cite{Ballerini2008,Ashraf2017,Lei2020,Mizumoto2025}. The coordination of the group then rests on the stability of each pair. In tandem running in blind termites, stability is defined in terms of maintaining tactile contact~\cite{Mizumoto2024,Mizumoto2025,Caroll2026}. The anatomy of termites, in turn, imposes spatiotemporal constraints on the physical contact required to preserve coordinated movement~\cite{Mizumoto2025}. Motivated by the empirical studies, we transform the dynamics of the pair into coordinates relative to the leader. In the new coordinates, translational degrees of freedom are eliminated, and the dynamics of the leader-follower pair reduce to radial and angular fluctuations on the complex plane.

\section*{Results}
\label{sec:main}
\subsection*{Dynamics of the herd}
\label{subsec:skew}
We first introduce the dynamics of the herding system on the state space and then transform to the coordinates relative to the leader.
\paragraph*{Leader}
The leader moves at a constant speed, denoted \(v_G\), and the direction of travel is subject to a Vicsek-type angular perturbation introduced by stochastic noise~\cite{Vicsek1995}. We replace the noise with the forcing produced by a chaotic map. The resulting motion remains irregular, while the system is deterministic.

We denote the position of the leader at time $t$ by $\boldsymbol{r}_G(t)\in\mathbb{R}^2$, with time step $\Delta t = 1$, and the direction of travel by $\theta(t)$.  The unit direction vector is
\(
\bm n\bigl(\theta(t)\bigr)
=[\cos\theta(t),\,\sin\theta(t)]^\top.
\)
Let \(F\) be a chaotic map. In the present work, we choose the Gauss~\cite{Wirsing1974}, logistic~\cite{May1976}, and standard~\cite{Chirikov1979} maps for \(F\). For the one-dimensional Gauss and logistic maps, \(F\colon[0,1]\to[0,1]\), and the state \(x(t)\) of the map lies in \([0,1]\) for every \(t\). For the standard map, let \(q(t),p(t)\in[0,1)\) be two coordinates of the state \(x(t)=(q(t),p(t))\). We define the projection \(h(q,p)=q\), so that \(h(x(t))=q(t)\in[0,1)\). We let \(\omega(t)\) denote the turning angle at time \(t\), bounded by a fixed parameter \(\omega_{\max}>0\). The value of \(F\) in the unit interval is mapped to \(\omega(t)\in[-\omega_{\max},\omega_{\max}]\). Finally, the dynamics of the leader are given by
\begin{equation}
\begin{aligned}
x(t) &= F\bigl(x(t-1)\bigr),\\
\omega(t) &= \omega_{\max}
\begin{cases}
2x(t)-1, & \text{Gauss, logistic},\\
2h(x(t))-1, & \text{standard},
\end{cases}\\
\theta(t) &= \theta(t-1)+\omega(t),\\
\bm r_G(t) &= \bm r_G(t-1)+v_G\bm n\bigl(\theta(t)\bigr).
\end{aligned}
\label{eq:generator}
\end{equation}
The definitions of the chaotic maps, state spaces, and invariant measures are given in
Materials and Methods.

\paragraph*{Predecessor selection}
The predecessor selection takes place in two stages. First, a follower considers up to a fixed number of nearest neighbors among the followers closer to the leader in Euclidean distance. Second, the follower follows the neighbor closest to the leader, or the leader directly when there is no such neighbor.
For \(i=1,\ldots,N\), let \(\bm r_R^i(t)\in\mathbb R^2\) denote the position of the \(i\)th follower. Define the Euclidean distance from the position of follower \(i\) at the previous time step to the current position of the leader by
\begin{equation}
d_i(t)
=
\bigl\|
\bm r_R^i(t-1)-\bm r_G(t)
\bigr\|.
\label{eq:distance}
\end{equation}
From \autoref{eq:distance}, the followers closer to the leader than follower \(i\) form the set
\(
\cEi=\{\,j:d_j(t)<d_i(t)\,\}.
\)
Let \(k\) be a fixed positive integer, and let \(|\cEi|\)
denote the number of followers in \(\cEi\). Define \(\cKik\subseteq\cEi\) as the set of the
\(\min\{k,|\cEi|\}\) followers nearest to follower \(i\) at time \(t-1\)~\cite{Ballerini2008}.
When \(\mathcal K_i^k(t)=\varnothing\), follower \(i\) follows
the leader directly. Otherwise, let \(j^*\) denote the index
of the follower in \(\mathcal K_i^k(t)\) closest to the leader,
\[
j^*
=
\arg\min_{j\in\mathcal K_i^k(t)}d_j(t).
\]
If distances are equal during either step of the selection, the follower with the smaller index is chosen.
The position of the predecessor is given by
\begin{equation}
\bm r_{\mathrm{prd}}^i(t)
=
\begin{cases}
\bm r_R^{j^*}(t-1),
& \mathcal K_i^k(t)\neq\varnothing,\\[4pt]
\bm r_G(t),
& \mathcal K_i^k(t)=\varnothing.
\end{cases}
\label{eq:predecessor}
\end{equation}
At each time step, the state of the system determines the predecessor tree rooted in the leader~\cite{BangJensen2009}.

\paragraph*{Subsystem of a follower as a convolution}
A follower is constructed as a subsystem, not just as a position vector. For \(i=1,\ldots,N\), follower \(i\) has an individual subsystem, which consists of three nonphysical states, the proxy \(\bm r_P^i(t)\), the auxiliary \(\bm r_A^i(t)\), and the transport \(\bm r_T^i(t)\). The three internal states provide a distributed delay, and the delay retains the influence of a previous predecessor despite the switch caused by \autoref{eq:predecessor}. 
Let \(\lambda,\mu,\delta\in(0,1)\) and \(\alpha,\gamma>0\) be fixed parameters.
The subsystem evolves sequentially within each time step, 
\begin{align}
\bm r_P^i(t)
&=
\bm r_P^i(t-1)
+\lambda\bigl(
\bm r_{\mathrm{prd}}^i(t)-\bm r_P^i(t-1)
\bigr),
\nonumber\\
\bm r_A^i(t)
&=
\mu\,\bm r_A^i(t-1)
+\gamma\bigl(
\bm r_P^i(t)-\bm r_R^i(t-1)
\bigr),
\nonumber\\
\bm r_T^i(t)
&=
\delta\,\bm r_T^i(t-1)
+\alpha\,\bm r_A^i(t),
\nonumber\\
\bm r_R^i(t)
&=
\bm r_R^i(t-1)+\bm r_T^i(t).
\label{eq:SubsysDym}
\end{align}
The sequential composition of the proxy, auxiliary, and transport maps is the PAT operator.
We define the state of the \(i\)th subsystem as
\(
\bm w^i(t)=
\bigl(\bm r_P^i(t),\bm r_A^i(t),
\bm r_T^i(t),\bm r_R^i(t)\bigr).
\)
Let \(A\in\mathbb R^{4\times4}\) be the state matrix and \(B\in\mathbb R^{4\times1}\) the coupling vector so that 
\begin{equation}
A=
\begin{pmatrix}
1-\lambda & 0 & 0 & 0\\
\gamma(1-\lambda) & \mu & 0 & -\gamma\\
\alpha\gamma(1-\lambda) & \alpha\mu & \delta & -\alpha\gamma\\
\alpha\gamma(1-\lambda) & \alpha\mu & \delta & 1-\alpha\gamma
\end{pmatrix},
\qquad
B=
\begin{pmatrix}
\lambda\\ \gamma\lambda\\ \alpha\gamma\lambda\\ \alpha\gamma\lambda
\end{pmatrix}.
\label{eq:AB}
\end{equation}
In matrix form, the subsystem in \autoref{eq:SubsysDym} follows the recursion
\begin{equation}
\bm w^i(t)
=
A\bm w^i(t-1)+B\bm r_{\mathrm{prd}}^i(t).
\label{eq:w}
\end{equation}
For the spectral radius \(\rho(A)\) of \(A\), the parameters are chosen so that \(\rho(A)<1\), as specified in
\autoref{eq:coupling}.

\paragraph*{Roles of the internal states.} Each of the five parameters has a specific role in forming the  following behavior. 
Let \(\lambda\in(0,1)\) denote the weight given to the current position of the predecessor in the proxy \(\bm r_P^i(t)\), so that \(1-\lambda\) sets how much of the past is retained. Over time, the initial value of the proxy fades out. The proxy approaches a weighted average of the positions of the current predecessor and of all earlier predecessors of follower \(i\). 

The auxiliary \(\bm r_A^i(t)\) is an exponentially weighted sum of current and past discrepancies. The vector \(\bm r_P^i(t)-\bm r_R^i(t-1)\) is the discrepancy between the updated proxy and the previous position of the follower. A larger \(\mu\) retains earlier discrepancies longer, whereas a larger \(\gamma\) increases the effect of the current discrepancy. The follower corrects the mismatch left by the previous correction, with an anticipation informed by the tendency of the motion of the predecessor.  

The transport \(\bm r_T^i(t)\) accumulates the auxiliary state in the same way, with \(\delta\) and \(\alpha\) in the roles analogous to those of \(\mu\) and \(\gamma\). The transport is one step the follower takes, and the size and direction of the step are set by the accumulated correction, not by the discrepancy at the moment.

The three states correspond to a perception, a processing, and an actuation: the proxy is what the follower has taken in of the predecessor, the auxiliary is what the follower makes of the gap, the transport is the motion that comes out. Through the PAT, the motion of the follower is neither purely reactive nor synchronized with the predecessor.

\subsubsection*{Change of coordinates}
The position of a follower relative to the leader has a radial and an angular component, and the complex plane represents the two as modulus and argument.
Let \(\mathrm{i}\) denote the imaginary unit, \(\mathrm{i}^2=-1\). We identify each planar vector \(\bm r=(r_x,r_y)\in\mathbb R^2\) with the complex number \(r_x+\mathrm{i}r_y\).
Let \(\mathbb S^1\) denote the unit circle in \(\mathbb C\). For a complex number \(w\), \(\overline w\) denotes the complex conjugate of \(w\). Define the rotation factor by
\begin{equation}
z(t)=e^{-\mathrm{i}\theta(t)}\in\mathbb{S}^1.
\label{eq:z}
\end{equation}
At time \(t\), the position of follower \(i\) in the frame of the leader is
\begin{equation}
D_R^i(t)=\bigl(\bm r_R^i(t)-\bm r_G(t)\bigr)z(t).
\label{eq:Newcordi}
\end{equation}
In the new coordinates, the leader is fixed at the origin, and the direction of its travel is aligned with the positive real axis. The factor \(z(t)\) in \autoref{eq:Newcordi} rotates each displacement \(\bm r_R^i(t)-\bm r_G(t)\) by \(-\theta(t)\).

\subsection*{The asymptotic spatial measure}
\label{subsec:measure}
We define the distribution of the positions \(D_R^i(t)\) by a measure. Fix \(T\geq1\). For any point \(D\in\mathbb C\), the probability measure that assigns unit mass to \(D\) is the \emph{Dirac measure} \(\delta_D\). The distribution of the recorded positions, in the form of a probability measure on \(\mathbb C\), assigns to a region \(B\subseteq\mathbb C\) the fraction of time the followers have occupied \(B\). Such a measure is called an \emph{occupation measure} in ergodic theory~\cite{Eckmann1985}. The finite-time occupation measure is
\begin{equation}
\nu_T =
\frac{1}{NT}
\sum_{i=1}^{N}
\sum_{t=0}^{T-1}
\delta_{D_R^i(t)}.
\label{eq:nuT}
\end{equation}
We say that \(\nu_T\) \emph{converges weakly}~\cite{Billingsley2013} to a measure \(\nu\) on \(\mathbb C\) as \(T\to\infty\) if, for every bounded continuous function \(f\colon\mathbb C\to\mathbb R\),
\begin{equation}
\int f\,d\nu_T \longrightarrow \int f\,d\nu .
\label{eq:weak}
\end{equation}
When the weak limit \(\nu\) exists, we call it the \emph{asymptotic spatial measure}. 

\begin{figure*}[t]
\centering
\includegraphics[width=\textwidth]{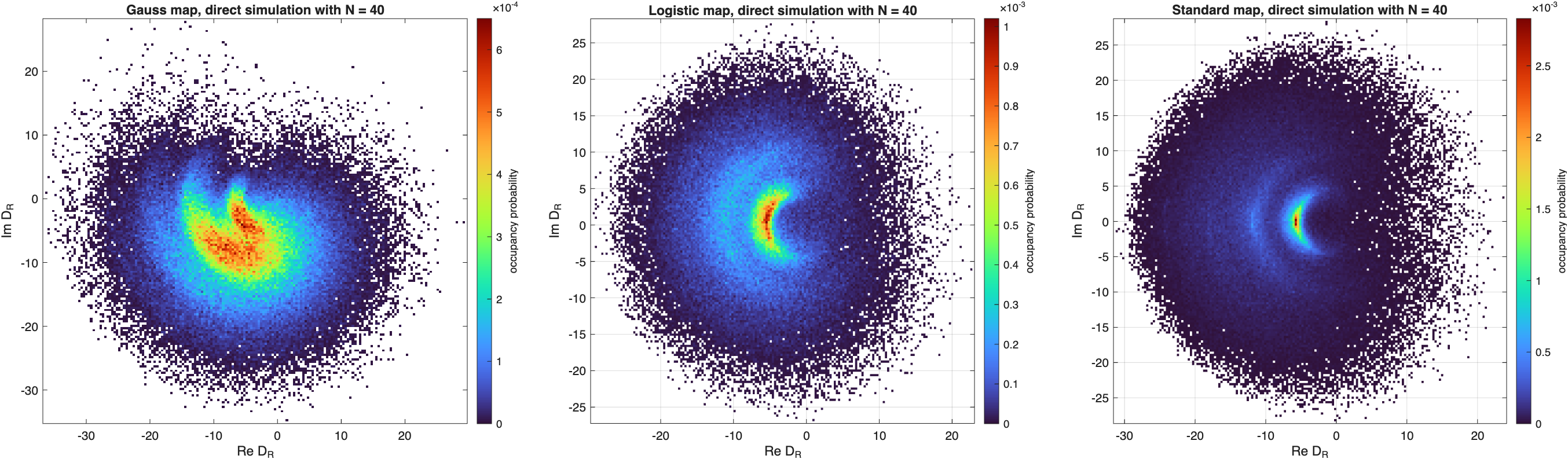}
\caption{Geometry of the finite-time occupation measure \(\nu_T\) from direct simulation with \(N=40\) followers and \(k=10\), for the Gauss, logistic, and standard generators. Each panel is computed from \(1.8\times10^{5}\) time steps of one recorded run. Each panel has its own bins, axis limits, and color scale.}
\label{fig:measure_simulation}
\end{figure*}

The position of the follower in the frame of the leader is a function of the history of the leader. The existence of \(\nu\) then follows from an invariant measure of \(F\). If \(F\) is also ergodic~\cite{Walters1982}, \(\nu\) is unique.

\(\nu\) is a probability measure on the plane, and its location, extent, and symmetry are the properties of a planar figure. The three together are called the geometric representation of \(\nu\). In the present paper, the geometric representation of \(\nu\) is the mathematical object we analyze. The representation of \(\nu_T\) for the three chaos generators is shown in \autoref{fig:measure_simulation} and indicates the dependence on the leader. The mass is symmetric about the direction of travel for the standard map, but lies on one side for the Gauss map. The logistic map falls between the two, with an asymmetry too weak to be seen. The asymmetry is discussed in the section on chirality.

We establish the existence and uniqueness of \(\nu\) for a single leader-follower pair, as shown in \autoref{si:measure}. The three chaotic maps we chose for \(F\), the Gauss, logistic, and standard maps, have known invariant measures, so \(\nu\) exists for each. Only the Gauss and logistic maps are ergodic~\cite{Walters1982,Lasota1994}, and \(\nu\) is unique for the two.\footnote{For the standard map at \(K=5\), regular islands exist~\cite{Duarte1994} and occupy about \(2\%\) of the phase space, where \(\nu\) differs (Materials and Methods).}

\subsection*{The fundamental relation and the selection operators}
\label{subsec:fundamental_kernel}

\paragraph*{The fundamental relation as an operator} The subsystem of each follower is linear in the position of the predecessor. In the leader-follower pair, the trajectory of the follower is then the image of the trajectory of the leader under a linear operator. We refer to the operator as the \emph{fundamental relation}. The distributed delay of the subsystem makes the fundamental relation a convolution over the positions of the predecessors up to the current time. Let \(\bm w^{(1)}\) denote the subsystem of the follower in the pair that defines the fundamental relation. The coordinate vector \(\bm e_R=(0,0,0,1)^\top\) leaves out the internal states from the subsystem, so that
\(
\bm r_R^{(1)}(t)=\bm e_R^\top\bm w^{(1)}(t)\). We denote the fundamental relation by \(\mathcal R_1\). For an integer \(\ell\geq0\), called the \emph{lag}, the influence of the predecessor from \(\ell\) time steps in the past enters the follower with a weight. We define the \emph{weight factor} by
\begin{equation}
\kappa_\ell=\bm e_R^\top A^\ell B.
\label{eq:position_kernel_main}
\end{equation}
The weight factor decays geometrically and the initial value of the subsystem fades out, both by \autoref{eq:coupling}.\footnote{The convolution series converge absolutely at each time \(t\). The weight factors decay geometrically by \autoref{eq:coupling}, and the leader moves at constant speed, so \(\|\bm r_G(t-\ell)\|\) grows at most linearly in \(\ell\).} Iterating the recursion in \autoref{eq:w} gives
\begin{equation}
\begin{aligned}
\bm r_R^{(1)}(t)
&=\bm e_R^\top\bm w^{(1)}(t)\\
&=\sum_{\ell=0}^{\infty}
\bm e_R^\top A^\ell B\,\bm r_G(t-\ell)\\
&=\sum_{\ell=0}^{\infty}
\kappa_\ell\bm r_G(t-\ell)\\
&=(\mathcal R_1\bm r_G)(t).
\end{aligned}
\label{eq:fundamental_response_operator}
\end{equation}
\autoref{eq:fundamental_response_operator} is the fundamental relation in explicit form. The operator is the convolution of the trajectory of the leader with the weight factors \(\kappa_\ell\).

\paragraph*{The fundamental relation in a system of \(N\) followers}
\label{subsec:nbody_operator}
We demonstrate that the dynamics of the herd with \(N\) followers are given by the operator of the pair composed with two new operators. The composition is the structure that makes an exact analysis possible. What the herd adds is the choice of predecessors, and the selection is where the complexity of the motion comes from.
By the interaction rules in \autoref{eq:predecessor}, the predecessor can be a different individual at each time step, either the leader or another follower. We specify the selection by two operators, one for each case. Index the leader by \(0\) and the followers by \(1,\ldots,N\). For follower \(i\), let \(p_i(t)\in\{0,1,\ldots,N\}\) denote the index of the predecessor selected at time \(t\). Let \(\bm r_R\) denote the collection of trajectories of all \(N\) followers. The indicator function \(\bm 1_{\{\cdot\}}\) is \(1\) when the enclosed condition is met and \(0\) otherwise. For a given sequence of predecessor indices \((p_i(t))_{i,t}\), define the \emph{leader-selection operator} \(\mathcal J\) and the \emph{follower-selection operator} \(\mathcal S\) by

\begin{equation}
\begin{aligned}
(\mathcal J\bm r_G)_i(t)
&=
\bm 1_{\{p_i(t)=0\}}\bm r_G(t),\\
(\mathcal S\bm r_R)_i(t)
&=
\begin{cases}
\bm r_R^{p_i(t)}(t-1), & p_i(t)\neq0,\\
0, & p_i(t)=0.
\end{cases}
\end{aligned}
\label{eq:selection_operators}
\end{equation}
The two operators correspond to the two cases of \autoref{eq:predecessor}, so that
\begin{equation}
\bm r_{\mathrm{prd}}^i(t)
=
(\mathcal J\bm r_G)_i(t)
+
(\mathcal S\bm r_R)_i(t).
\label{eq:prd_selection}
\end{equation}
The operator \(\mathcal R_1\) acts simultaneously on the sequence of predecessors of each follower, which gives
\begin{equation}
\bm r_R
=
\mathcal R_1
\bigl(
\mathcal J\bm r_G
+
\mathcal S\bm r_R
\bigr).
\label{eq:nbody_kernel}
\end{equation}
\autoref{eq:nbody_kernel} verifies that the fundamental relation is composed with the selection operators.

\subsection*{Exact factorization of the \texorpdfstring{$N$}{N}-body dynamics}
\label{subsec:factorization}

\paragraph*{Factorization with remainder}
\autoref{eq:nbody_kernel} has an exact consequence, the factorization of the operator. For \(N>1\), the dynamics written as an operator take the form of a sum of products of the fundamental relation and the selection operators. We refer to the sum as the \emph{factorization}. Let \(Q\in\mathbb N_0=\{0,1,2,\ldots\}\) denote the number of times \autoref{eq:nbody_kernel} is substituted into itself. Then, we obtain
\begin{align}
\bm r_R
&=
\mathcal R_1\mathcal J\bm r_G
+
\mathcal R_1\mathcal S
\mathcal R_1\mathcal J\bm r_G
+\cdots \nonumber\\
&\quad+
(\mathcal R_1\mathcal S)^Q
\mathcal R_1\mathcal J\bm r_G
+
(\mathcal R_1\mathcal S)^{Q+1}\bm r_R
\nonumber\\
&=
\sum_{q=0}^{Q}
(\mathcal R_1\mathcal S)^q
\mathcal R_1\mathcal J\bm r_G
+
(\mathcal R_1\mathcal S)^{Q+1}\bm r_R
\nonumber\\
&=
\mathcal R_1\Biggl[
\sum_{q=0}^{Q}
(\mathcal S\mathcal R_1)^q
\mathcal J\bm r_G
+
(\mathcal S\mathcal R_1)^Q
\mathcal S\bm r_R
\Biggr].
\label{eq:absolute_path_factorization}
\end{align}

We call \(Q\) the \emph{order} of the factorization, and \((\mathcal R_1\mathcal S)^{Q+1}\bm r_R\) the \emph{remainder} at order \(Q\). In a static predecessor tree, \(\mathcal R_1\) and \(\mathcal S\) commute, and \(\mathcal S^N=0\) by \autoref{eq:selection_operators} and \autoref{eq:predecessor}. The factorization then terminates by \(Q=N-1\). In the present system, the trees form an adaptive network so that it is not commutative, \(\mathcal R_1\mathcal S\neq\mathcal S\mathcal R_1\), in general. As a result, the terms with more than \(N-1\) selections are possible; see \autoref{si:fixedtree}.

\begin{figure}[t]
\centering
\includegraphics[width=\columnwidth]{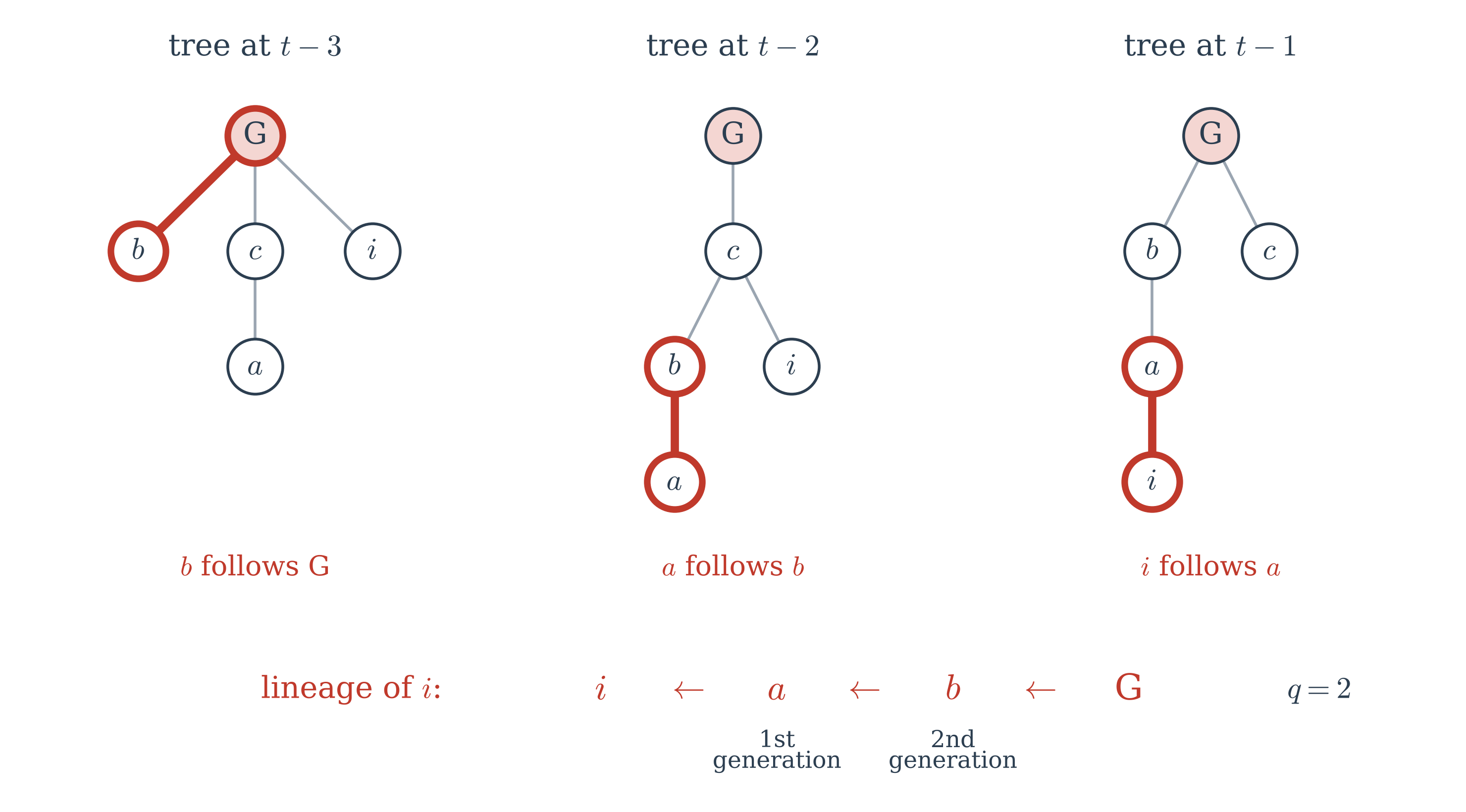}
\caption{A lineage of follower \(i\) across the predecessor trees at three consecutive times, for an illustrative sequence of predecessor indices. Each tree has the leader G at the root. At each time, the lineage takes one link (red). The predecessors \(a\) and \(b\) between \(i\) and the leader are the two generations of the lineage, \(q=2\). For simplicity, each step back along the lineage is one time step.}
\label{fig:tree}
\end{figure}

\begin{figure*}[t]
\centering
\includegraphics[width=0.98\textwidth]{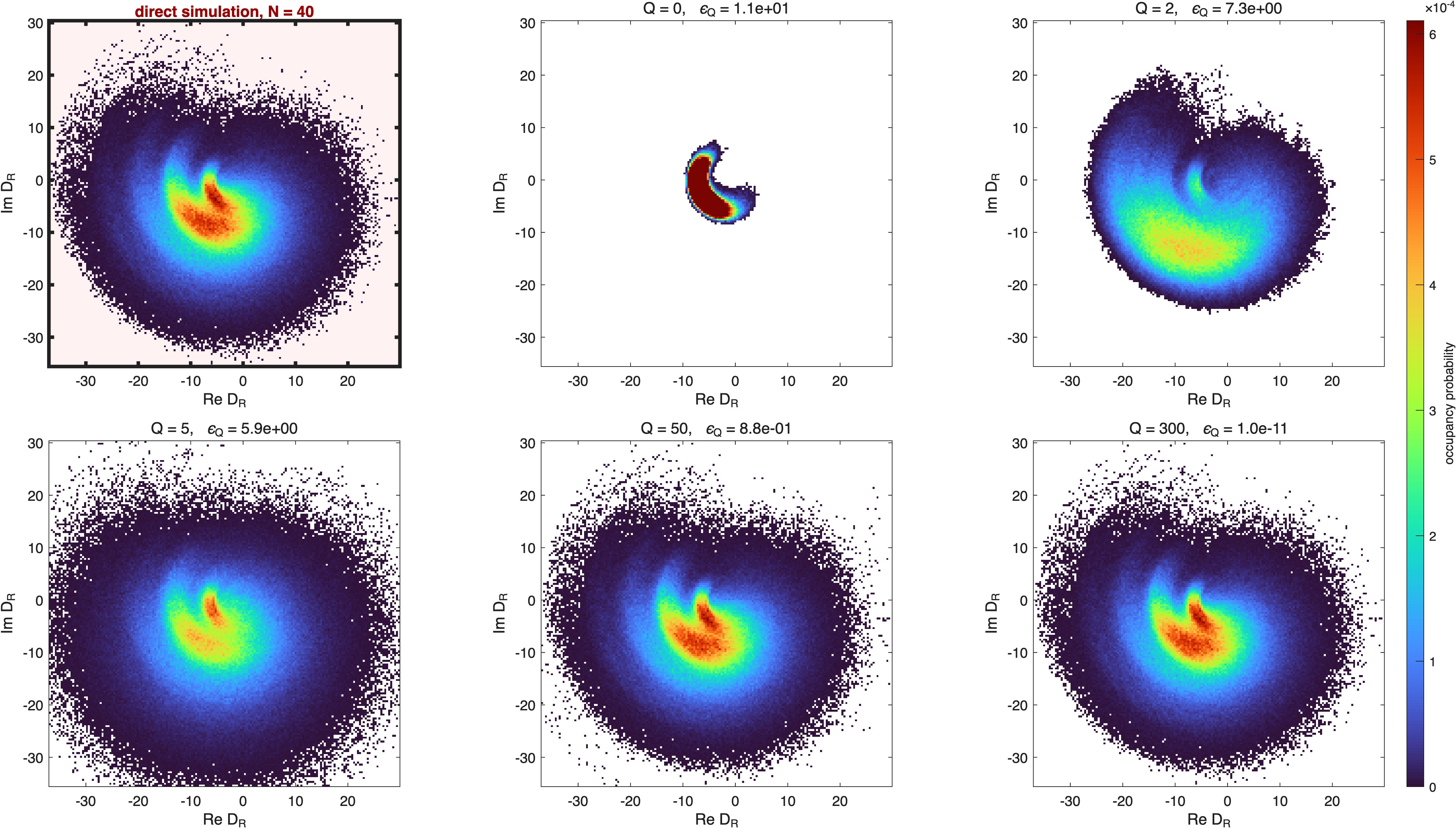}
\caption{Geometry of the finite-time occupation measure for the Gauss generator with $N=40$ followers and $k=10$. Top left: direct simulation. Remaining panels: factorization from $\mathcal R_1$ in \autoref{eq:factor} at order $Q$, with the remainder $\epsilon_Q$ of \autoref{eq:remainder_convergence}. All panels use the bins, axis limits, and color scale of the simulation, without fitting, translation, rotation, or rescaling. The factorization and the remainder are evaluated on $1.8\times10^{5}$ time steps of the recorded run.}
\label{fig:gauss_path_order}
\end{figure*}

\begin{figure*}[t]
\centering
\includegraphics[width=\textwidth]{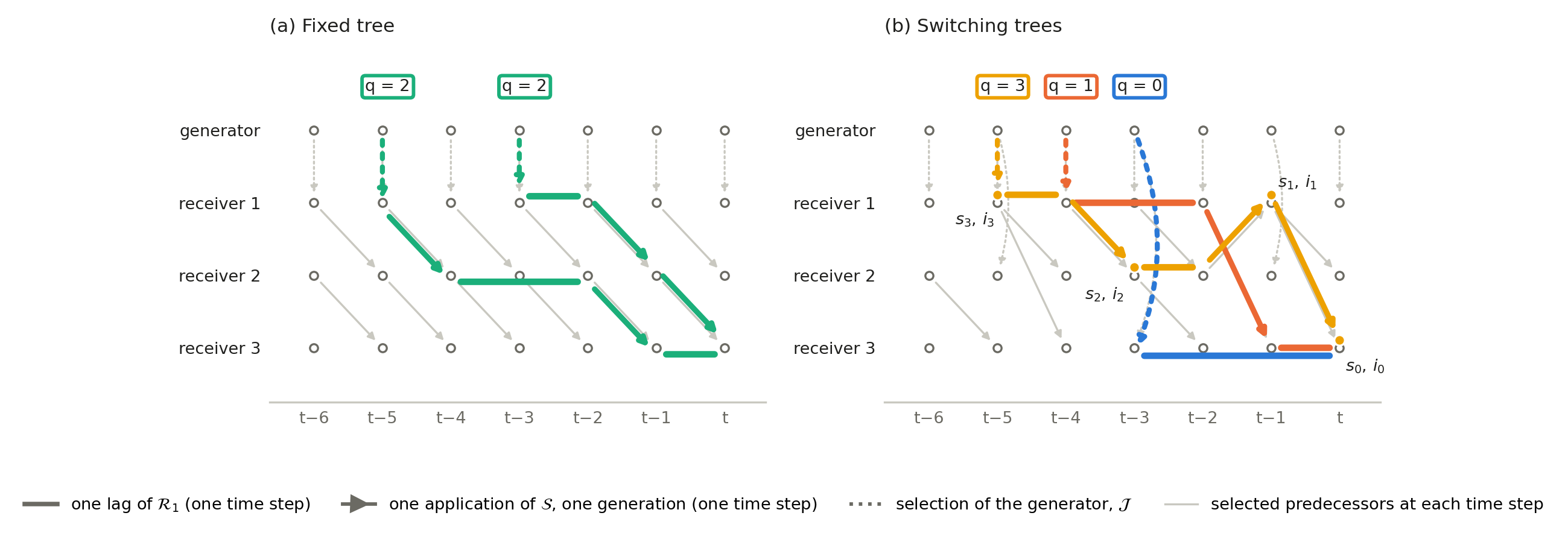}
\caption{Lineages to follower 3 at time \(t\) for \(N=3\), with predecessors chosen for illustration. Gray arrows mark the predecessor selected by each follower at each time step. A colored path is one lineage, read from the leader to follower 3, and its color marks the number of generations \(q\). (a) Fixed tree. Every lineage to follower 3 has \(q=2=N-1\). (b) Switching trees. Lineages to follower 3 have \(q=0\), \(1\), and \(3\). Along the lineage with \(q=3\), labeled by \(s_m\) and \(i_m\), follower 1 appears twice, follower 1 precedes follower 2 at \(t-3\), and follower 2 precedes follower 1 at \(t-1\).}
\label{fig:lineage}
\end{figure*}

\paragraph*{Expansion in relative coordinates and convergence}

Let \(\bm y=(y_i(t))_{i,t}\) denote any collection of complex numbers, one for each follower \(i\) and time \(t\). The rotation factor \(z\) in \autoref{eq:z} acts on \(\bm y\) at each time, \((z\bm y)_i(t)=z(t)y_i(t)\).
The displacement of the followers from the leader is

\(
\bm d=
\bigl(\bm r_R^i(t)-\bm r_G(t)\bigr)_{i,t},
\)
so that, by \autoref{eq:Newcordi},
\[
\bm D_R=z\bm d
=
\bigl(D_R^i(t)\bigr)_{i,t}.
\]
Assigning the trajectory of the leader to every follower defines a collection, \(y_i(t)=\bm r_G(t)\) for all \(i\) and \(t\), or simply \(\bm r_G\). We extend \(\mathcal J\) to collections by \((\mathcal J\bm y)_i(t)=\bm 1_{\{p_i(t)=0\}}y_i(t)\). Let \(I\) be the identity operator. Then, \autoref{eq:absolute_path_factorization} gives, for every finite \(Q\),
\begin{equation}
\bm D_R
=
z\sum_{q=0}^{Q}
(\mathcal R_1\mathcal S)^q
\bigl[\mathcal R_1(\mathcal J+\mathcal S)-I\bigr]\bm r_G
+
z(\mathcal R_1\mathcal S)^{Q+1}\bm d. 
\label{eq:factor}
\end{equation}

Let \(\|\cdot\|_{N,T}\) denote the root-mean-square norm over the \(N\)
followers and \(T\) time steps. Under the remainder condition, 
\begin{equation}
\lim_{Q\to\infty}
\left\|(\mathcal R_1\mathcal S)^{Q+1}\bm d\right\|_{N,T}=0,
\label{eq:remainder_convergence}
\end{equation}
the sum in \autoref{eq:factor} converges to \(\bm D_R\) with respect
to \(\|\cdot\|_{N,T}\). The factorization is exact for arbitrary \(N\) provided that the remainder decays to zero. The convolution operator
\(\mathcal R_1\) is independent of the size of the system, while the evolving rooted
predecessor tree specifies the combinatorial structure dependent on \(N\).

\subsection*{Numerical verification}
\label{subsec:numero}

\paragraph*{Comparison with direct simulation}
Numerical verification compares \(\bm D_R\) reconstructed from the
factorization in \autoref{eq:factor} with the one obtained from direct simulation for the system with 40 followers.

For the present work, the operators \(\mathcal J\), \(\mathcal S\), and
\(\mathcal R_1\) are defined analytically, while the trajectory of the leader \(\bm r_G\), rotation factor \(z\), and indices of the predecessors \(p_i(t)\) are obtained from direct simulation. The recorded indices specify the time-dependent action of \(\mathcal J\) and \(\mathcal S\), and the collection
\(\bigl(p_i(t)\bigr)_{i,t}\) defines the sequence of predecessor indices. The comparison tests the factorization for the recorded sequence of predecessor indices. 
\paragraph*{Residual and remainder}
Across the three generators, the shared parameters of the system are
\(
k=10,  v_G=\omega_{\max}=0.5,
 (\lambda,\mu,\delta,\gamma,\alpha)
=(0.7,0.85,0.8,0.025,0.1).
\)
Each simulation runs \(T=10^7\) iterates after a transient of
\(T_{\mathrm{tr}}=10^6\). The values of the parameters for the logistic and standard maps are \(r=4\) and \(K=5\), respectively. At order \(Q\) of the factorization, let \(\epsilon_{\rm id}\) denote
the root-mean-square norm of the difference between the left and
right-hand sides of \autoref{eq:factor}, including the remainder. Let \(\epsilon_Q\)
denote the norm of the finite-\(Q\) remainder appearing in
\autoref{eq:remainder_convergence}. \(\epsilon_{\mathrm{id}}\) checks the accuracy of the computation of the factorization, and \(\epsilon_Q\) is the quantity whose decay \autoref{tab:path_order} reports.

\begin{table}[t]
\caption{Numerical residual \(\epsilon_{\mathrm{id}}\) and finite-order
remainder \(\epsilon_Q\) for the \(N=40\) factorization in \autoref{eq:factor}.}

\centering

\begin{tabular}{lcccc}
Generator & \(Q\) & \(\epsilon_{\mathrm{id}}\) & \(\epsilon_Q\) & \(Q_{\mathrm{mach}}\) \\
\midrule
Gauss & 300 & \(3.612\times10^{-7}\) & \(1.022\times10^{-11}\) & 390 \\
Logistic & 50 & \(3.090\times10^{-7}\) & \(5.996\times10^{-10}\) & 78 \\
Standard & 50 & \(3.520\times10^{-7}\) & \(1.621\times10^{-16}\) & 50 \\
\bottomrule
\end{tabular}

\label{tab:path_order}
\end{table}
The values of $Q$ in \autoref{tab:path_order} are the orders at which
$\epsilon_Q$ falls below the corresponding numerical residual
$\epsilon_{\mathrm{id}}$ for each generator.
As \(Q\) increases further, \(\epsilon_Q\) continues to decay and reaches double-precision machine epsilon at \(Q=390\), \(78\), and \(50\) for the Gauss, logistic, and standard generators, respectively. All three orders exceed \(N-1=39\), the largest number of selections within a single tree. The complete \(Q\)-dependence in
\autoref{fig:convergence} shows the decay of
\(\epsilon_Q\) with increasing \(Q\). The corresponding
geometric comparisons appear in \autoref{fig:gauss_path_order},
\autoref{fig:logistic_path_order_supp}, and
\autoref{fig:standard_path_order_supp}.  
Our result takes the sequence of predecessors as given, recorded from the simulation. For a given sequence, the trajectory of the leader alone determines the trajectories of all followers through the fundamental relation. A rigorous characterization of the sequence of predecessor indices and their combinatorial effects on the geometry of the system lies beyond the scope of the current paper.

\paragraph*{Complexity of propagation}

The order of the factorization quantifies how far the dependence on the leader reaches into the herd along the sequence. The order counts the followers through which the influence of the leader passes before reaching a follower. Fix a follower \(i\) and a time \(t\). Counting backward in time from \(i\) to the leader forms a sequence of individuals. We call the sequence a \emph{lineage} of \(i\), and each individual between \(i\) and the leader a \emph{generation} of the lineage. The number of generations is denoted by \(q\). A lineage with two generations is shown in \autoref{fig:tree}. The term of index \(q\) in \autoref{eq:absolute_path_factorization} collects the lineages of \(q\) generations.

In a static tree, the lineage of every follower has at most \(N-1\) generations, one for each selection. \(Q\) exceeds the bound in \autoref{tab:path_order}. Lineages in a static tree and trees of the adaptive network are shown in \autoref{fig:lineage}. In \autoref{eq:absolute_path_factorization}, a lineage of \(q\) generations corresponds to one product, in which the fundamental relation acts \(q+1\) times and the follower selection \(q\) times. The two alternate,
\[
(\mathcal R_1\mathcal S)^q\mathcal R_1\mathcal J\bm r_G
=\underbrace{\mathcal R_1\mathcal S\,\mathcal R_1\mathcal S\cdots\mathcal S\,\mathcal R_1}_{q+1\text{ actions of }\mathcal R_1,\ q\text{ of }\mathcal S}\,\mathcal J\bm r_G .
\]

Each action of the fundamental relation is a temporal convolution, which combines all past positions of the selected predecessors with the weight factors \(\kappa_\ell\) in \autoref{eq:position_kernel_main}. The influence along a lineage of \(q\) generations has an explicit form in terms of the past positions of the leader. Let \(\bm\ell=(\ell_0,\ldots,\ell_q)\) denote the lags of the \(q+1\) actions. 
Let \(s_m\) denote the selection times and \(i_m\) the followers along the lineage,
\begin{align*}
s_0&=t-\ell_0, & s_{m+1}&=s_m-1-\ell_{m+1},\\
i_0&=i, & i_{m+1}&=p_{i_m}(s_m).
\end{align*}

For a given choice of lags \(\bm\ell\), the lineage from \(i\) reaches the leader after some number of generations, not necessarily \(q\). Let \(\chi_i^{(q)}\) denote the indicator that the lineage reaches the leader after exactly \(q\) generations,
\[
\chi_i^{(q)}(t;\bm\ell)=\prod_{m=0}^{q-1}\bm 1_{\{p_{i_m}(s_m)\neq0\}}\;\bm 1_{\{p_{i_q}(s_q)=0\}}.
\]
Then
\begin{equation}
\begin{aligned}
\bigl[(\mathcal R_1\mathcal S)^q\mathcal R_1\mathcal J\bm r_G\bigr]_i(t)
&=\sum_{\ell_0,\ldots,\ell_q\geq0}\kappa_{\ell_0}\cdots\kappa_{\ell_q}\,
\chi_i^{(q)}(t;\bm\ell)\\
&\quad\times\bm r_G\Bigl(t-\sum_{m=0}^{q}\ell_m-q\Bigr).
\end{aligned}
\label{eq:lineage_term}
\end{equation}
A past position of the leader reaches the follower through every combination of \(q+1\) lags, weighted by the product \(\kappa_{\ell_0}\cdots\kappa_{\ell_q}\), so that the weights are in \((q+1)\)-linear form. The tree of the adaptive network extends the propagation beyond \(N-1\) generations. The influence of one follower on another can become locally bidirectional over time, even though the global coupling rule is strictly one-way. Neither a single trajectory nor the network at a single instant can show the property.

\subsection*{Chirality of the geometry}
\label{subsec:chirality}

The factorization shows that, for a given sequence of predecessor indices, the trajectory of the leader determines the geometric representation of the measure. The influence of the leader shows up in another property of the representation, the handedness. We provide the explanation by analyzing the symmetry of the equations of motion and of the map. 

An object is chiral if it cannot be superimposed on its mirror image. Suppose every turn of the leader is reversed, \(\omega\to-\omega\). Suppose also that all initial positions are mirrored across the initial direction of travel. The predecessor selection rule in \autoref{eq:predecessor} depends only on distances, which are preserved under reflection. 

The PAT recursion in \autoref{eq:w} is likewise linear with real scalar coefficients and does not distinguish left from right\footnote{Reflection across the direction of travel is a linear map \(M\) on the plane. Since \(A\) and \(B\) are real and act on the two components alike, applying \(M\) to \autoref{eq:w} gives \(M\bm w^i(t)=A\,M\bm w^i(t-1)+B\,M\bm r_{\mathrm{prd}}^i(t)\), so the mirrored states obey the same recursion.} Any handedness of the motion with respect to the direction of travel must come from the leader~\cite{Curie1894,Chalmers1970}. 

The standard map commutes with \((q,p)\mapsto(1-q,1-p)\), which reverses \(\omega\). The leader with the standard generator is achiral. The invariant measure of the Gauss map is asymmetric. The Gauss generator produces a chiral leader already in the mean turn, \(\omega_{\max}(2/\log 2-3)\approx-0.115\,\omega_{\max}\). 

Each turn of the logistic generator, taken alone, is achiral, since the distribution of \(\omega(t)\) is symmetric. However, the dependence between consecutive turns produces chirality. Consecutive turns satisfy \(\omega(t+1)/\omega_{\max}=1-2[\omega(t)/\omega_{\max}]^2\). The average of \(\omega(t)^2\omega(t+1)\) over the invariant measure is then \(-\omega_{\max}^3/4\). The average is odd under the reversal \(\omega\to-\omega\) and vanishes for an achiral leader~\cite{Brillinger1965,Nikias1993}. The nonzero value indicates that the leader with the logistic generator is chiral. The next turn after a large turn is on average toward negative \(\omega\). 

For the pair, we compute the imaginary part of the mean of \(\nu\), normalized by the standard deviation of \(\nu\) (\autoref{si:chirality}). The value is \(2\times10^{-4}\) for the standard generator, \(-0.544\) for the Gauss generator, and \(0.081\) for the logistic generator. The influence of the leader is thus observable in the location and extent of the measure, in the reach of the factorization, and in the handedness of the measure.

\section*{Discussion}
\label{sec:conclude}

For the leader-follower pair that forms the fundamental relation, the asymptotic spatial measure exists for almost every initial condition of the leader (\autoref{si:measure}). The result shows that a statistical measure of nonequilibrium dynamics can be obtained from chaotic forcing, without introducing external noise. The only source of irregularity is \(F\), which enters through the turning angle of the leader. Yet, despite the irregularity and the repeated switching of the predecessors, each follower moves smoothly and the herd does not break apart. The memory effect of PAT, the contraction property of the matrix \(A\), and the two-step predecessor selection rule might contribute to coherence of the herd.

The asymptotic statistics of the herd are computable, as the only source of irregularity is the map \(F\). The three chaotic maps selected for \(F\) have known invariant measures, so that the statistics follow from ergodic theory. The mean and the variance of \(\nu\) are computed from the invariant measure of the map, or as Birkhoff averages along one orbit when the invariant measure is not in closed form.  There are two notable properties in our construction that are relevant to the computability. The irregularity enters through a chaotic driver that is unaffected by the herd, and the memory of each member fades. Our approach may extend to other collective motions, and what the two properties produce here may inform constructions with further properties.

The existence of the measure is not proven for the herd with more than one follower. In simulations, the occupation measure continues to change throughout the run. We have observed that the change is slow compared with the memory of a follower. For example, the memory of a follower lasts about \(22\) time steps by the contraction \(\rho(A)\approx0.954\), and the relative position of the pair depends on the history of the leader with exponentially fading memory.

For the logistic generator, the variance of the measure decreases from about \(158\) to \(104\) over the first \(4\times10^6\) time steps. The variance remains near \(104\) up to \(1.3\times10^7\) time steps. At that point, the orbit of the logistic map in double precision reaches the fixed point \(x=0\), and the leader turns by \(-\omega_{\max}\) at every time step. Beyond \(1.3\times10^7\) time steps, the simulation no longer represents the chaotic forcing, and the convergence is inconclusive. For the Gauss generator, the variance decreases from about \(109\) to \(83\) over \(2\times10^7\) time steps, with an abrupt change near \(1.65\times10^7\). The convergence is inconclusive over the period. One explanation might be that the slow time scale arises from the relative positions of the followers, which set the trees that can form and change only through rare reorganizations. The existence of the measure for the herd and the time scale of the reorganizations are left for future work.

In a moving animal group, chirality is often attributed to an anatomical property of the individual. In flying geese and schooling squid, the position of an individual relative to others has been associated with a preference for one eye~\cite{Zaynagutdinova2022,Sakurai2024}. Tandem running of blind termites, however, relies only on contact~\cite{Mizumoto2021}, so an eye preference could not explain chirality. In our system, the follower has no lateral bias, yet for the logistic and Gauss generators the mean position of the follower is displaced to one side, and the measure is chiral. We show that the displacement comes from the turns of the leader, through their mean for the Gauss map and through the dependence between consecutive turns for the logistic map (\emph{Chirality of the geometry}). The chirality of a moving group could well be dynamical in origin. The hypothesis can be tested with the available trajectory data of tandem runs~\cite{Mizumoto2021,Mizumoto2025}.

The factorization concerns all followers at once. For a given sequence of predecessor indices, the positions of all followers follow exactly from the fundamental relation, and the structure of the network enters only through the lineages in \autoref{eq:lineage_term}. The current work provides an adaptive network with exactly solvable dynamics. The existence of one such network is the result; a general means of analysis may develop from understanding the case.

In the current work, the sequence of predecessor indices is selected from the positions of the followers, and the positions of the followers are computed from the sequence. A sequence that the dynamics can produce must select itself. Closing the factorized dynamics leads to a fixed-point problem for the sequence of predecessor indices. An unsolved inverse problem is to recover all the sequences of predecessor indices that produce a given measure. Discovering a hidden structure depends not only on analyzing the system in detail but also on choosing the right object to observe.

\matmethods{

\subsection*{Generator maps}
\label{subsec:mm_generators}
The state space \(X\) is specified for each map. For the Gauss generator, \(X=[0,1]\) and
\begin{equation*}
F_{\mathrm G}(s)=
\begin{cases}
\dfrac{1}{s}-\left\lfloor\dfrac{1}{s}\right\rfloor, & s>0,\\
0, & s=0.
\end{cases}
\end{equation*}
The Gauss map preserves the probability measure~\cite{Wirsing1974}
\(d\mathbb P_{\mathrm G}(s)=ds/[(1+s)\log 2]\). Since the mean of \(s\) with respect to \(\mathbb P_{\mathrm G}\) is not \(1/2\), the corresponding turning angles have nonzero mean
\(\mathbb E_{\mathbb P_{\mathrm G}}[\omega]=\omega_{\max}(2/\log 2-3)\neq0\).
For the logistic generator, \(X=[0,1]\) and \(F_{\mathrm L}(s)=4s(1-s)\). The invariant probability measure is \(d\mathbb P_{\mathrm L}(s)=ds/[\pi\sqrt{s(1-s)}]\). Its symmetry about \(s=1/2\) gives \(\mathbb E_{\mathbb P_{\mathrm L}}[\omega]=0\). Equivalently, under the change of variables \(s=\sin^2(\pi y)\), the logistic map is semiconjugate to the doubling map \(y\mapsto 2y\pmod 1\), and \(\omega(y)=-\omega_{\max}\cos(2\pi y)\).
For the standard generator, \(X=\mathbb T^2\) and \(x=(q,p)\). In normalized coordinates, the map is
\begin{equation*}
\begin{aligned}
p(t)&=p(t-1)+\frac{K}{2\pi}\sin\bigl(2\pi q(t-1)\bigr)\pmod 1,\\
q(t)&=q(t-1)+p(t)\pmod 1,
\end{aligned}
\end{equation*}
with \(K=5\). The standard map preserves Lebesgue measure on \(\mathbb T^2\). For Lebesgue measure, \(\mathbb E_{\mathbb P_{\mathrm S}}[\omega]=0\).

\subsection*{Stability of the follower recursion}
\label{subsec:mm_stability}
The recursion in \autoref{eq:SubsysDym} has the matrix form \(\bm w^i(t)=A\bm w^i(t-1)+B\bm r_{\mathrm{prd}}^i(t)\), with \(A\) and \(B\) in \autoref{eq:AB}.
The first row of \(A\) has a single nonzero entry, so \(A\) is block lower triangular. The eigenvalues of \(A\) are \(1-\lambda\) and the roots of the characteristic polynomial of the lower \(3\times3\) block. The Jury conditions for a cubic polynomial~\cite{Jury1964} show that \(\rho(A)<1\) if and only if
\begin{equation}
\alpha\gamma<\min\left\{2(1+\mu)(1+\delta),\ \frac{(1-\mu)(1-\delta)(1-\mu\delta)}{\mu\delta}\right\}
\label{eq:coupling}
\end{equation}
(\autoref{si:stability}). For our parameters, \(\alpha\gamma=0.0025\), below the bound \(0.0141\). Under \autoref{eq:coupling} there exist \(C>0\) and \(\eta\in(\rho(A),1)\) with \(\|A^t\|\leq C\eta^t\)~\cite{Horn1985}, and the contribution of the initial follower state decays exponentially for any predecessor history.

\subsection*{Existence of the spatial measure}
\label{subsec:mm_existence}
For the leader-follower pair that forms the fundamental relation, define the tail weights
\begin{equation}
\tau_\ell=\sum_{\ell'>\ell}\kappa_{\ell'}=\bm e_R^\top A^{\ell+1}(I-A)^{-1}B,\qquad \ell\geq0.
\label{eq:tau}
\end{equation}
Since \(\sum_\ell\kappa_\ell=1\) and \(\bm r_G(t)-\bm r_G(t-1)=v_G\overline{z(t)}\), summation by parts gives
\begin{equation}
D_R(t)=-v_G\sum_{\ell=0}^{\infty}\tau_\ell\,z(t)\overline{z(t-\ell)}+\varepsilon_t,\qquad|\varepsilon_t|\leq C_0\eta^t,
\label{eq:Dsum}
\end{equation}
for constants \(C_0>0\) and \(0<\eta<1\). The term \(\varepsilon_t\) is the contribution from the initial follower state and from the leader positions before \(t=0\). Let \((\widehat X,\widehat{\mathbb P},\widehat F)\) be the natural extension of \((X,\mathbb P,F)\), where \(\widehat x=(\ldots,x_{-1},x_0)\) satisfies \(F(x_j)=x_{j+1}\), and let \(\psi\) denote the rescaling in \autoref{eq:generator}, so that \(\theta(t)-\theta(t-1)=\psi(x_t)\). Define the relative-coordinate observable
\begin{equation}
\mathcal D(\widehat x)=-v_G\sum_{\ell=0}^{\infty}\tau_\ell\exp\Bigl(-\mathrm{i}\sum_{j=-\ell+1}^{0}\psi(x_j)\Bigr),
\label{eq:Dobs}
\end{equation}
so that \(D_R(t)=\mathcal D(\widehat F^t\widehat x)+\varepsilon_t\). Absolute summability of \((\tau_\ell)\) bounds \(|\mathcal D|\), so all relative positions lie in a compact disk. The pointwise ergodic theorem on the natural extension~\cite{Walters1982}, applied to a countable dense family of continuous test functions, gives convergence of the time averages for \(\widehat{\mathbb P}\)-almost every history, and the Riesz representation theorem~\cite{Rudin1987} identifies the limit with a probability measure \(\nu\). If \(F\) is ergodic, \(\nu=\mathcal D_{\#}\widehat{\mathbb P}\). The full proof is given in \autoref{si:measure}.

\subsection*{Operator equations}
\label{subsec:mm_operators}
Throughout the current section, we consider histories for which the lag sums converge absolutely and the convolution representation holds without a homogeneous initial-state term. The follower-position component is \(\bm r_R^i(t)=\sum_{\ell\geq0}\kappa_\ell\bm r_{\mathrm{prd}}^i(t-\ell)\). Separating the two possible predecessor types gives
\[
\begin{aligned}
\bm r_R^i(t)
&=\sum_{\ell\geq0}\kappa_\ell\bm 1_{\{p_i(t-\ell)=0\}}\bm r_G(t-\ell)\\
&\quad+\sum_{\ell\geq0}\kappa_\ell\bm 1_{\{p_i(t-\ell)\neq0\}}\bm r_R^{p_i(t-\ell)}(t-\ell-1).
\end{aligned}
\]
The first sum is \((\mathcal R_1\mathcal J\bm r_G)_i(t)\), and the second is \((\mathcal R_1\mathcal S\bm r_R)_i(t)\). Hence, \autoref{eq:nbody_kernel} follows directly for the given sequence of predecessor indices, and \autoref{eq:absolute_path_factorization} follows by induction on \(Q\). With \(\bm d=\bm r_R-\bm r_G\),
\[
\bm d
=\mathcal R_1\mathcal J\bm r_G+\mathcal R_1\mathcal S(\bm d+\bm r_G)-\bm r_G
=\bigl[\mathcal R_1(\mathcal J+\mathcal S)-I\bigr]\bm r_G+\mathcal R_1\mathcal S\bm d,
\]
and repeated substitution followed by pointwise multiplication by \(z(t)\) gives \autoref{eq:factor}. Fix a Schur-stable PAT recursion shared by all followers, and suppose that every follower has one predecessor at each time. For arbitrary \(N\), any leader trajectory, and any sequence of predecessor indices, \autoref{eq:factor} holds at every finite \(Q\). Consequently, the finite factorization is independent of the leader dynamics, number of followers, neighborhood structure, and predecessor-selection rule. Under \autoref{eq:remainder_convergence},
\begin{equation}
\bm D_R=z\sum_{q=0}^{\infty}(\mathcal R_1\mathcal S)^q\bigl[\mathcal R_1(\mathcal J+\mathcal S)-I\bigr]\bm r_G.
\label{eq:complete}
\end{equation}

\subsection*{Numerical evaluation}
\label{subsec:mm_numerical}
Let \(\bm D_R^{\mathrm{dir}}\) denote the relative-coordinate histories obtained by direct simulation. In terms of the norm \(\|\cdot\|_{N,T}\),
\[
\epsilon_{\mathrm{id}}=\bigl\|\bm D_R^{\mathrm{dir}}-\bm D_{R,\mathrm{id}}^{(Q)}\bigr\|_{N,T},
\qquad
\epsilon_Q=\bigl\|z(\mathcal R_1\mathcal S)^{Q+1}\bm d\bigr\|_{N,T},
\]
where \(\bm D_{R,\mathrm{id}}^{(Q)}\) is the complete right-hand side of \autoref{eq:factor}, including its remainder. The operators \(\mathcal R_1\) and \((\mathcal R_1\mathcal S)^{Q+1}\) are applied over the full recorded history, and the norms are evaluated on a window of \(2\times10^{5}\) time steps after a burn-in of \(2\times10^{4}\) steps.
\begin{figure}[t]
\centering
\includegraphics[width=0.95\linewidth]{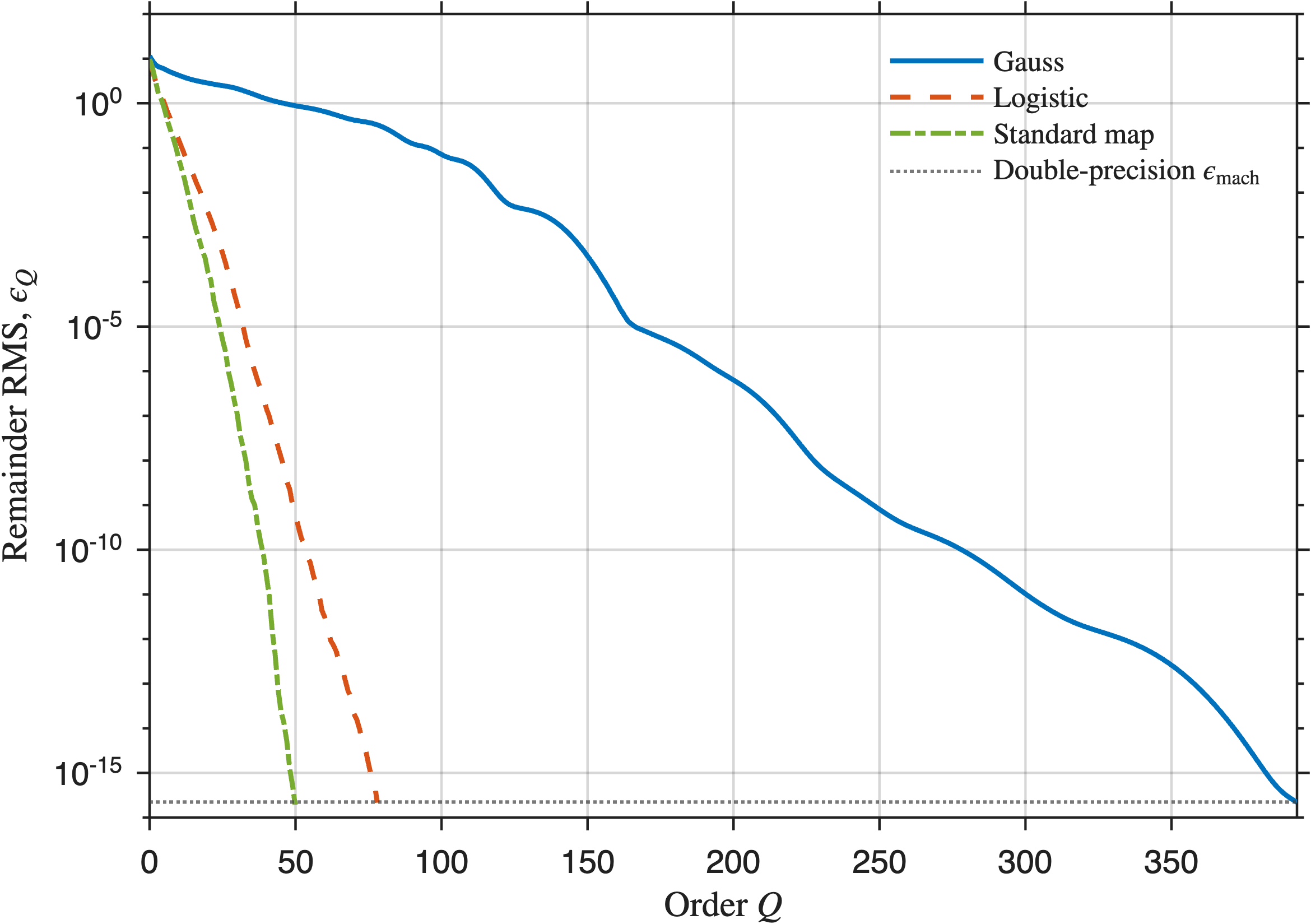}
\caption{Decay of the factorization remainder for the Gauss, logistic, and standard-map generators with \(N=40\). At each order \(Q\), \(\epsilon_Q\) is the root-mean-square norm defined in Materials and Methods. The dotted line marks double-precision machine epsilon. The curves numerically evaluate \autoref{eq:remainder_convergence} for the recorded sequence of predecessor indices; the finite-\(Q\) identity does not depend on the limit.}
\label{fig:convergence}
\end{figure}

\subsection*{Asymptotic behavior of the occupation measure}
\label{subsec:mm_longtime}
The long-time behavior is examined in an independent implementation of the model in Python, with the parameters of the main text and forty followers. The leader starts at the origin with \(\theta(0)=0\), and the followers start at positions drawn from a normal distribution with standard deviation \(2\) around the origin, with proxy states equal to the positions and auxiliary and transport states equal to zero. After a transient of \(2\times10^5\) time steps, each run continues for \(2\times10^7\) time steps, and the mean and variance of the relative positions are computed in blocks of \(5\times10^5\) time steps. Over \(5\times10^6\) time steps, four runs from different initial conditions give occupation measures that agree within the sampling error for the Gauss generator, and for the standard-map generator when the initial condition lies in the chaotic region. For the standard map at \(K=5\), about \(2\%\) of initial conditions lie in regular islands, and such a run gives a different measure, consistent with the dependence of \(\nu\) on the history when \(F\) is not ergodic (\autoref{si:measure}).

\subsection*{AI assistance}
\label{subsec:mm_ai}
Claude (Anthropic) assisted with editing, LaTeX preparation, coding, and algebraic and numerical checks. The author takes full responsibility for the content and writing.

}
\showmatmethods{}

\dataavail{Simulation code and the recorded sequence of predecessor indices will be deposited in a public repository upon publication.}

\acknow{The author thanks Ira B. Schwartz for posing the research question that initiated the study and for many helpful discussions. She thanks her collaborator Nobuaki Mizumoto, whose work on tandem running in termites has been indispensable to the study. She thanks Martin Golubitsky for suggesting the simplest-case analysis that shaped the present result. She thanks Maurizio Porfiri for contributing to the initial form of the model and Claude (Anthropic) for assistance with editing, LaTeX preparation, coding, and algebraic and numerical checks. Lastly, she thanks James Yorke for working through part of the manuscript line by line with her, not to make the paper intelligible but to make her able to write intelligible text, and to recognize when it is not, with the help of AI. The initial development of the present work was supported by the Office of Naval Research under Grant No. N0001425GI01182 and No. N0001425GI01158.}
\showacknow{}

\bibliography{Kamimoto_preprint}

\clearpage
\onecolumn
\appendix
\begin{center}{\Large\sffamily\bfseries Supplementary Information}\end{center}
\setcounter{equation}{0}\renewcommand{\theequation}{S\arabic{equation}}
\setcounter{table}{0}\renewcommand{\thetable}{S\arabic{table}}
\setcounter{figure}{0}\renewcommand{\thefigure}{S\arabic{figure}}
\sisection{Stability of the follower recursion}
\label{si:stability}

 Write the per-follower state as \(w^i(t)=(r_P^i(t),r_A^i(t),r_T^i(t),r_R^i(t))\in\mathbb R^4\) (one-dimensional notation; the two-dimensional case is identical componentwise). For each follower \(i\), each \(t\geq1\), and any predecessor history \(r_{\mathrm{prd}}^i(\cdot)\), the follower state satisfies \(w^i(t)=Aw^i(t-1)+Br_{\mathrm{prd}}^i(t)\), where \(A\) and \(B\) are given in \autoref{eq:AB}.

\subsection*{Eigenvalue structure}
The first row of \(A\) has a single nonzero entry, so the proxy row decouples from the rest, and \(A\) is block lower triangular. Its eigenvalue is \(1-\lambda\), which satisfies \(|1-\lambda|<1\) for all \(\lambda\in(0,1)\). The remaining \(3\times 3\) block governing \((r_A^i,r_T^i,r_R^i)\) is
\begin{equation*}
M=\begin{pmatrix}
\mu & 0 & -\gamma\\
\alpha\mu & \delta & -\alpha\gamma\\
\alpha\mu & \delta & 1-\alpha\gamma
\end{pmatrix}.
\end{equation*}
The eigenvalues of \(A\) are \(1-\lambda\) and the eigenvalues of \(M\). The characteristic polynomial of \(M\) is
\begin{equation}
p(\xi)=\xi^3-(1+\mu+\delta-\alpha\gamma)\xi^2+(\mu+\delta+\mu\delta)\xi-\mu\delta.
\label{si:charpoly}
\end{equation}

\subsection*{Schur stability conditions}
All roots of \(p(\xi)\) lie strictly inside the unit disk if and only if the following four conditions hold. They are the Jury (Schur-Cohn) conditions for a cubic polynomial~\cite{Jury1964}. Write \(p(\xi)=\xi^3+a_2\xi^2+a_1\xi+a_0\), where
\[
a_2=-(1+\mu+\delta-\alpha\gamma),\qquad
a_1=\mu+\delta+\mu\delta,\qquad
a_0=-\mu\delta.
\]
\textbf{C1:} \(p(1)>0\). Computing directly, \(p(1)=\alpha\gamma>0\), satisfied for all \(\alpha,\gamma>0\).

\textbf{C2:} \(p(-1)<0\). Computing directly, \(p(-1)=-(2+2\mu+2\delta+2\mu\delta-\alpha\gamma)\). Thus C2 holds if and only if \(\alpha\gamma<2(1+\mu)(1+\delta)\).

\textbf{C3:} \(|a_0|<1\) where \(a_0=-\mu\delta\). This gives \(\mu\delta<1\), satisfied for all \(\mu,\delta\in(0,1)\).

\textbf{C4:} \(1-a_0^2>|a_0a_2-a_1|\). Since \(\mu\delta<1\), one computes \(|a_0a_2-a_1|=(\mu+\delta)(1-\mu\delta)+\alpha\gamma\mu\delta\). Since \(1-a_0^2=(1-\mu\delta)(1+\mu\delta)\), condition C4 is equivalent to \((1-\mu)(1-\delta)(1-\mu\delta)>\alpha\gamma\mu\delta\).

Conditions C1 and C3 are satisfied automatically for \(\lambda,\mu,\delta\in(0,1)\) and \(\alpha,\gamma>0\). Conditions C2 and C4 require the bounded coupling condition, \autoref{eq:coupling}. Under this condition, all eigenvalues of \(A\) lie strictly inside the unit disk, and hence \(\rho(A)<1\).

\subsection*{Decay of the initial state}
For each follower \(i\), each \(t\geq1\), and any predecessor history \(r_{\mathrm{prd}}^i(\cdot)\), the follower state satisfies
\begin{equation}
w^i(t)=A^tw^i(0)+\sum_{s=1}^{t}A^{t-s}B\,r_{\mathrm{prd}}^i(s).
\label{si:variation}
\end{equation}
Moreover, there exist constants \(C>0\) and \(0<\eta<1\), depending only on \(A\), such that
\begin{equation}
\|A^t\|\leq C\eta^t\quad\text{for all }t\geq0,
\label{si:decay}
\end{equation}
hence \(\|A^tw^i(0)\|\leq C\eta^t\|w^i(0)\|\to0\). Consequently,
\begin{equation}
w^i(t)-\sum_{s=1}^{t}A^{t-s}B\,r_{\mathrm{prd}}^i(s)\to0,\qquad t\to\infty.
\label{si:fading}
\end{equation}

\paragraph*{Proof}
Iterating \(w^i(t)=Aw^i(t-1)+Br_{\mathrm{prd}}^i(t)\) from the initial condition \(w^i(0)\) gives \(w^i(1)=Aw^i(0)+Br_{\mathrm{prd}}^i(1)\) and \(w^i(2)=A^2w^i(0)+ABr_{\mathrm{prd}}^i(1)+Br_{\mathrm{prd}}^i(2)\). By induction, \(w^i(t)=A^tw^i(0)+\sum_{s=1}^{t}A^{t-s}Br_{\mathrm{prd}}^i(s)\). By Gelfand's formula~\cite{Horn1985}, since \(\rho(A)<1\), there exist constants \(C>0\) and \(\eta\in(\rho(A),1)\) such that \(\|A^t\|\leq C\eta^t\) for all \(t\geq0\), proving \autoref{si:decay}. Applying \autoref{si:decay} to the homogeneous term gives \(\|A^tw^i(0)\|\leq\|A^t\|\,\|w^i(0)\|\leq C\eta^t\|w^i(0)\|\to0\). Rearranging \autoref{si:variation} gives \(w^i(t)-\sum_{s=1}^{t}A^{t-s}Br_{\mathrm{prd}}^i(s)=A^tw^i(0)\to0\), which is \autoref{si:fading}. The decay rate depends only on \(A\), hence only on the structural parameters \(\lambda,\mu,\delta,\alpha,\gamma\).

\sisection{Existence of the asymptotic spatial measure}
\label{si:measure}
We show the existence result stated in the main text. For the leader-follower pair that forms the fundamental relation, \(N=1\), and the finite-time occupation measure of \autoref{eq:nuT} reduces to \(\nu_T=T^{-1}\sum_{t=0}^{T-1}\delta_{D_R(t)}\). Weak convergence is defined in \autoref{eq:weak}.

Let \(F:X\to X\) denote the generator map and let \(\mathbb P\) be an invariant Borel probability measure for \(F\). Let \((\widehat X,\widehat{\mathbb P},\widehat F)\) be its natural extension, where \(\widehat x=(\ldots,x_{-1},x_0)\) satisfies \(F(x_j)=x_{j+1}\). The shift is \(\widehat F(\ldots,x_{-1},x_0)=(\ldots,x_0,F(x_0))\). The turning angle is determined by the bounded measurable function \(\psi\), the rescaling in \autoref{eq:generator}, so that \(\theta(t)-\theta(t-1)=\psi(x_t)\), \(x_t=x(t)\).

For the pair, write \(w(t)=Aw(t-1)+Br_G(t)\) and \(r_R(t)=e_R^\top w(t)\), where \(e_R=(0,0,0,1)^\top\), and assume \(\rho(A)<1\). Define the follower coefficients \(\kappa_{\ell}=e_R^\top A^{\ell}B\), \({\ell}\geq0\). Stability of \(A\) implies exponential decay of \(\kappa_{\ell}\). The follower coefficients satisfy \(\sum_{\ell\geq0}\kappa_{\ell}=e_R^\top(I-A)^{-1}B=1\). The tail weights \(\tau_\ell\) are defined in \autoref{eq:tau}. Exponential decay of \(A^{\ell}\) implies \(\sum_{\ell=0}^{\infty}|\tau_\ell|<\infty\).

Iteration of the stable recursion of the pair gives
\[
r_R(t)=\sum_{\ell=0}^{\infty}\kappa_{\ell}r_G(t-\ell)+\widetilde\varepsilon_t,\qquad|\widetilde\varepsilon_t|\leq C_0\eta^t,
\]
for some \(C_0>0\) and \(0<\eta<1\). The term \(\widetilde\varepsilon_t\) is the contribution from the initial follower state and from the leader positions before \(t=0\) that the infinite sum includes. It decays exponentially, because \(\kappa_\ell\) decays exponentially while \(r_G\) grows at most linearly.

Since \(z(t)=e^{-\mathrm{i}\theta(t)}\), its complex conjugate is \(\overline{z(t)}=e^{\mathrm{i}\theta(t)}\). The leader displacement satisfies \(r_G(t)-r_G(t-1)=v_G\overline{z(t)}\). With \(\sum_\ell\kappa_\ell=1\), the relative position is \(r_R(t)-r_G(t)=\sum_{\ell\geq0}\kappa_\ell[r_G(t-\ell)-r_G(t)]+\widetilde\varepsilon_t\). Writing each difference as a sum of one-step displacements and exchanging the order of summation gives \(-\sum_{\ell\geq0}\tau_\ell[r_G(t-\ell)-r_G(t-\ell-1)]\). Rearranging the absolutely convergent sums and multiplying by \(z(t)\) gives \autoref{eq:Dsum}. The rotation factors satisfy
\[
z(t)\overline{z(t-\ell)}=\exp\Bigl(-\mathrm{i}\sum_{j=t-\ell+1}^{t}\psi(x_j)\Bigr).
\]
The sum in the exponent is zero for \(\ell=0\). With the observable \(\mathcal D\) of \autoref{eq:Dobs}, \autoref{eq:Dsum} reads \(D_R(t)=\mathcal D(\widehat F^t\widehat x)+\varepsilon_t\).

Absolute summability of \((\tau_\ell)\) gives \(|\mathcal D(\widehat x)|\leq M\) with \(M=v_G\sum_{\ell=0}^{\infty}|\tau_\ell|\). Both \(\mathcal D(\widehat F^t\widehat x)\) and \(D_R(t)\) belong to the compact disk \(K=\{D\in\mathbb C:|D|\leq M+C_0\}\).

Compactness guarantees that the occupation measures have weakly convergent subsequences. Convergence of the full sequence follows from the pointwise ergodic theorem applied to the natural extension of the leader dynamics~\cite{Walters1982}. Choose a countable uniformly dense family \(\{\varphi_m\}_{m\geq1}\subset C(K)\), where \(C(K)\) is the space of continuous real functions on \(K\). For every \(m\), the pointwise ergodic theorem applied to \(\varphi_m\circ\mathcal D\) gives convergence of \(T^{-1}\sum_{t=0}^{T-1}\varphi_m(\mathcal D(\widehat F^t\widehat x))\) for \(\widehat{\mathbb P}\)-almost every \(\widehat x\). Intersecting the corresponding full-measure sets gives simultaneous convergence for every \(m\). Uniform density extends the limit to every \(\varphi\in C(K)\).

Define
\[
L_{\widehat x}(\varphi)=\lim_{T\to\infty}\frac1T\sum_{t=0}^{T-1}\varphi\bigl(\mathcal D(\widehat F^t\widehat x)\bigr).
\]
The functional \(L_{\widehat x}\) is linear and positive, with \(L_{\widehat x}(1)=1\). The Riesz representation theorem~\cite{Rudin1987} gives a unique Borel probability measure \(\nu_{\widehat x}\) on \(K\) satisfying \(L_{\widehat x}(\varphi)=\int_K\varphi(D)\,d\nu_{\widehat x}(D)\). Since \(\varepsilon_t\to0\), uniform continuity of \(\varphi\) on \(K\) gives
\[
\lim_{T\to\infty}\frac1T\sum_{t=0}^{T-1}\bigl|\varphi(D_R(t))-\varphi\bigl(\mathcal D(\widehat F^t\widehat x)\bigr)\bigr|=0.
\]
Consequently,
\begin{equation}
\nu_T\rightharpoonup\nu_{\widehat x}\qquad\text{as }T\to\infty
\label{si:weak}
\end{equation}
for \(\widehat{\mathbb P}\)-almost every leader history \(\widehat x\). Hence the convergence required in the definition of weak convergence holds, and \(\nu_{\widehat x}\) is the asymptotic spatial measure of the pair.

If \(F\) is ergodic with respect to \(\mathbb P\), the limiting measure is independent of the typical leader history and satisfies \(\nu=\mathcal D_{\#}\widehat{\mathbb P}\), \(\nu(B)=\widehat{\mathbb P}(\mathcal D^{-1}(B))\) for every Borel set \(B\subseteq\mathbb C\). Without ergodicity, the limiting measure may depend on the ergodic component containing the leader history, while weak convergence still holds for \(\widehat{\mathbb P}\)-almost every history.

\sisection{Chirality of the pair}
\label{si:chirality}
For the pair, \autoref{eq:Dsum} gives the mean of the asymptotic spatial measure \(\nu\),
\begin{equation}
m_{(1)}
=
\int_{\mathbb C} D\,d\nu(D)
=
-v_G\sum_{\ell=0}^{\infty}\tau_\ell\,\mathbb E\bigl[z(t)\overline{z(t-\ell)}\bigr],
\label{eq:m1_supp}
\end{equation}
where \(\mathbb E\) denotes the average over the histories of the leader with respect to \(\widehat{\mathbb P}\). Let \(v_{(0)}=\int_{\mathbb C}|D-m_{(1)}|^2\,d\nu(D)\) denote the variance of \(\nu\). When every turn of the leader is reversed, \(\omega\to-\omega\), each average \(\mathbb E[z(t)\overline{z(t-\ell)}]\) is replaced by its complex conjugate, and so is \(m_{(1)}\). For an achiral generator, every such average is real, and \(\operatorname{Im}m_{(1)}=0\). A nonzero \(\operatorname{Im}m_{(1)}\) shows that \(\nu\) is chiral.

\autoref{tab:chirality_supp} lists \(m_{(1)}\), \(v_{(0)}\), and \(\operatorname{Im}m_{(1)}/\sqrt{v_{(0)}}\) for the three generators with the parameters of the main text. For the Gauss map, the averages in \autoref{eq:m1_supp} are evaluated with a Chebyshev collocation of the Gauss-Kuzmin-Wirsing transfer operator on \(80\) nodes, with the infinite sum over inverse branches represented analytically beyond \(1024\) explicit branches. For the logistic map at \(r=4\), the conjugacy \(x=\sin^2(\pi y)\) turns the map into \(y\mapsto 2y \bmod 1\) with Lebesgue measure, and the averages are evaluated by propagating Fourier coefficients under the doubling map with Bessel-function multipliers, with cutoff \(32\). Both evaluations are converged to \(10^{-10}\). For the standard map, the averages are evaluated as weighted Birkhoff averages over \(1.1\times10^{7}\) recorded iterates of the \(N=1\) run at \(K=5\). The lag sum is taken to \(\ell=560\), where the omitted tail is below \(4\times10^{-11}\). For the logistic map, the average at lag one is real, and the first nonzero imaginary part appears at lag two, in agreement with the origin of the chirality in consecutive turns.

\begin{table}[t]
\caption{Mean \(m_{(1)}\), variance \(v_{(0)}\), and normalized imaginary mean of the asymptotic spatial measure of the pair. The standard-map value is a Birkhoff average, and its imaginary part lies within the sampling error.}
\centering
\begin{tabular}{lccc}
\toprule
Generator & \(m_{(1)}\) & \(v_{(0)}\) & \(\operatorname{Im}m_{(1)}/\sqrt{v_{(0)}}\) \\
\midrule
Standard & \(-4.979159+0.000692\,\mathrm i\) & \(10.3779\) & \(2\times10^{-4}\) \\
Logistic & \(-4.556573+0.302642\,\mathrm i\) & \(13.9795\) & \(0.081\) \\
Gauss    & \(-5.798338-1.619016\,\mathrm i\) & \(8.8519\)  & \(-0.544\) \\
\bottomrule
\end{tabular}
\label{tab:chirality_supp}
\end{table}

\sisection{Fixed-tree termination}
\label{si:fixedtree}
Suppose the predecessor tree is independent of time. Since every follower predecessor is closer to the leader, the tree is acyclic and rooted at the leader. Repeated application of \(\mathcal S\) follows successive follower predecessors, and the action is zero once the leader is reached. Since there are \(N\) followers, \(\mathcal S^N=0\).

For a time-independent predecessor tree, \(\mathcal S\) acts by a fixed follower-selection matrix together with a one-step time shift. The operator \(\mathcal R_1\) is the same temporal convolution on every follower component, and so commutes with \(\mathcal S\). Hence \((\mathcal R_1\mathcal S)^N=\mathcal R_1^N\mathcal S^N=0\). Consequently, the remainder \((\mathcal R_1\mathcal S)^{Q+1}\bm d\) of \autoref{eq:factor} vanishes for \(Q\geq N-1\). Lineages of more than \(N-1\) generations can occur only when predecessor selection changes in time.

\sisection{Paths along a lineage}
\label{si:paths}
Let \(T=\sum_{m=0}^{q}\ell_m+q\) denote the delay between follower \(i\) at time \(t\) and the past position \(\bm r_G(t-T)\) of the leader in \autoref{eq:lineage_term}. Going back from time \(t\) to time \(t-T\), each unit step either adds one to a lag \(\ell_m\) or passes to the predecessor through \(\mathcal S\). Exactly \(q\) of the \(T\) steps are selections. Each choice of \(q\) selections among the \(T\) steps determines exactly one vector of lags \(\bm\ell\). For order \(q\), the count is \(\binom{T}{q}\), and the sum over all orders is \(2^T\). With the selection times \(s_m\) and the followers \(i_m\) along the lineage, as defined before \autoref{eq:lineage_term}, the indicator
\begin{equation}
\chi_i^{(q)}(t;\bm\ell)=\prod_{m=0}^{q-1}\bm 1_{\{p_{i_m}(s_m)\neq0\}}\;\bm 1_{\{p_{i_q}(s_q)=0\}}
\label{eq:indicator_supp}
\end{equation}
equals one only for a lineage that reaches the leader after exactly \(q\) selections, and zero otherwise. At order \(q\), no more than \(\binom{T}{q}\) terms of follower \(i\) contain \(\bm r_G(t-T)\).

Let \(M=B\bm e_R^\top\) denote the matrix that passes the output of a predecessor to the input of its follower. Each weight \(\kappa_\ell=\bm e_R^\top A^\ell B\) is a scalar, and the product of the weights regroups as
\begin{equation}
\kappa_{\ell_0}\kappa_{\ell_1}\cdots\kappa_{\ell_q}
=\bm e_R^\top A^{\ell_0}MA^{\ell_1}M\cdots MA^{\ell_q}B .
\label{eq:product_supp}
\end{equation}
The regrouping follows from the associativity of the matrix product. Between \(\bm e_R^\top\) and \(B\), every lag step contributes \(A\), and each selection \(M\), in temporal order. \(A\) represents the memory of PAT, and \(M\) the passage along the tree. Through \autoref{eq:indicator_supp}, the sequence of predecessor indices decides which ordered products of \(A\) and \(M\) appear in \autoref{eq:lineage_term}. In a static tree, a lineage contains at most \(N-1\) selections, so no product has more than \(N-1\) factors \(M\). Switching removes the bound, and a lineage may return to a follower it has left.

\sisection{Factorization for the logistic and standard generators}
\label{si:comparison}
The factorization at the orders of \autoref{tab:path_order} is compared with the direct simulation for the logistic and standard generators in \autoref{fig:logistic_path_order_supp} and \autoref{fig:standard_path_order_supp}.

\begin{figure*}[t]
\centering
\addtocounter{figure}{1}
\begin{minipage}[t]{0.03\textwidth}\sublabel{fig:logistic_path_order_supp}\end{minipage}%
\begin{minipage}[t]{0.92\textwidth}\vspace{0pt}\includegraphics[width=\textwidth]{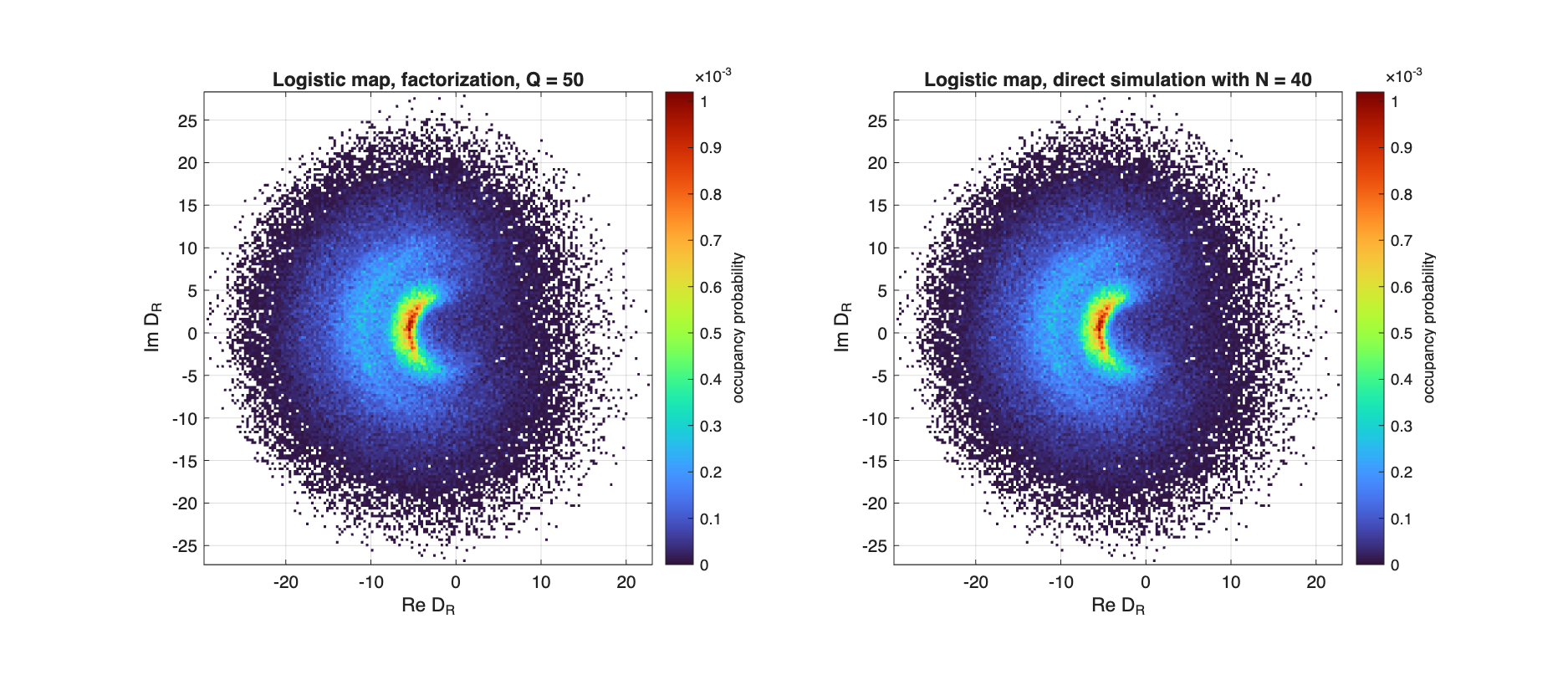}\end{minipage}\par\vspace{0.5em}
\begin{minipage}[t]{0.03\textwidth}\sublabel{fig:standard_path_order_supp}\end{minipage}%
\begin{minipage}[t]{0.92\textwidth}\vspace{0pt}\centering\includegraphics[width=0.85\textwidth]{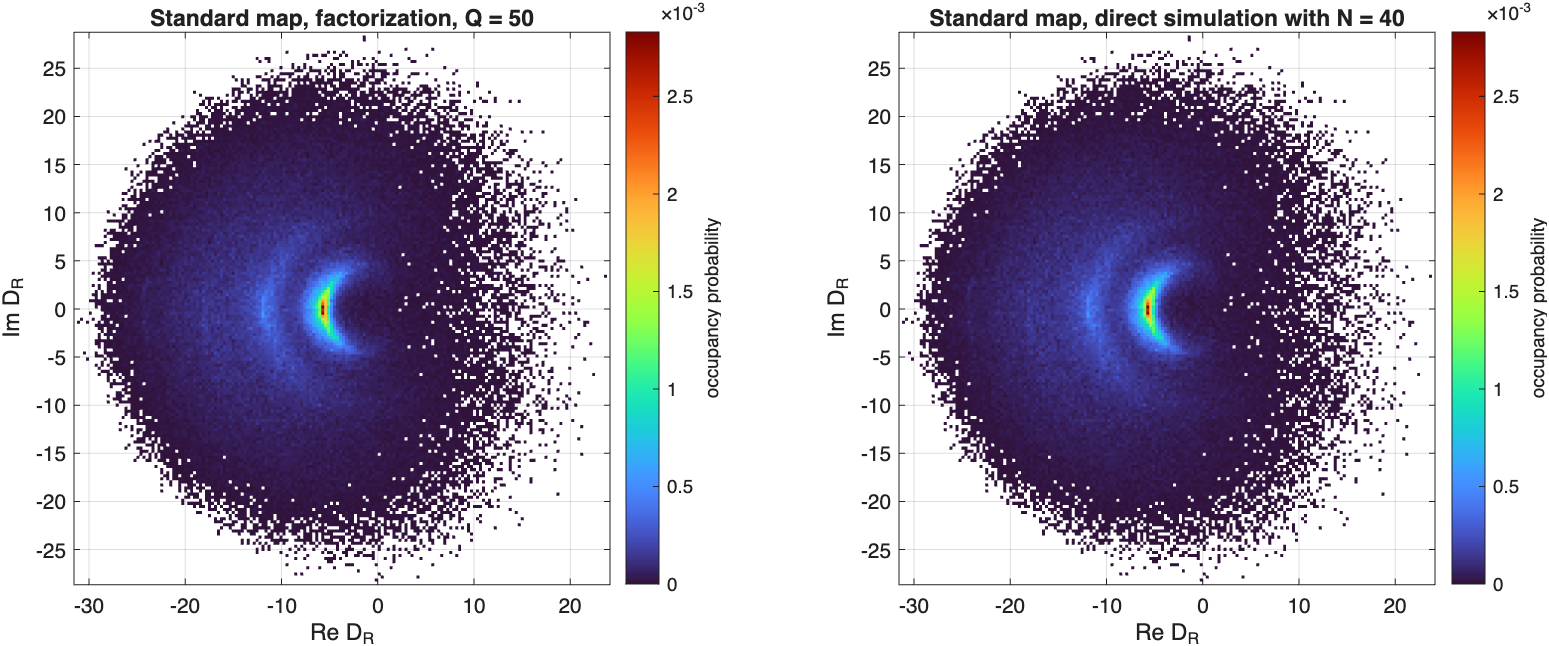}\end{minipage}
\addtocounter{figure}{-1}
\caption{Numerical evaluation of the factorization for the logistic and standard-map generators with \(N=40\) followers and \(k=10\). (a) Logistic generator, \(Q=50\). (b) Standard-map generator, \(Q=50\). In each row, the left panel is the factorization from \(\mathcal R_1\) in \autoref{eq:factor} at the order \(Q\) at which \(\epsilon_Q\) falls below \(\epsilon_{\mathrm{id}}\), and the right panel is the direct simulation. The two panels in each row share the bins, limits, and color scale, without fitting, translation, rotation, or rescaling.}
\label{fig:logistic_standard_path_order_supp}
\end{figure*}

\end{document}